\documentclass[11pt,reqno]{article}
\usepackage{cite}
\usepackage{amsmath,amsthm, amssymb}
\usepackage{amsfonts}
\numberwithin{equation}{section}
\title{\bf On  a representation of Humbert's double hypergeometric series $\Phi_3$ in a series of Gauss's $_2F_1$ function}  
\author{\bf \small  Arjun K. Rathie$^1$,  Victor V. Manako$^{2*}$, Harsh Vardhan Harsh $^3$ \\\\
{\footnotesize ${}^1$ Department of Mathematics, School of Mathematical and Physical Sciences,}\\ {\footnotesize Central University of Kerala,  Kasaragod- 671316,}\\ {\footnotesize Kerala - INDIA} \\
{\footnotesize E-Mail: akrathie@cukerala.edu.in} \\\\
{ \footnotesize ${}^2$ Department of General Physics and Physics of Oil and Gas Industry,}  \\ {\footnotesize Samara State Technical University, }\\
{ \footnotesize Molodogvardeyskaya Street,  244,} \\ {\footnotesize Samara, 443100, Russia} \\
{ \footnotesize E-Mail: victor.manako@mail.ru} \\\\
{ \footnotesize $^3$ Department of Mathematics,} \\ { \footnotesize Amity School of Engineering and Technology,} \\  { \footnotesize Amity University Rajasthan, NH-11C,} \\ { \footnotesize Kant Kalwar, Jaipur-303002, India.}} 
\begin{document}
\date{}
\maketitle
\begin{abstract}
Very recently a new series representation of Humbert's double hypergeometric series $\Phi_3$ in series of Gauss's $_2F_1$ function was  given by one of us. The aim of this short research note is to provide an alternative proof of the result. A few interesting special cases are also given.\vspace{0.4cm}

\noindent {\bf Mathematics Subject Classification:}  33C15, 33E20, 35C70

\vspace{0.3cm}
\noindent {\bf Keywords:} Hypergeometric function; Humbert's function.
\end{abstract}
\vspace{0.3cm}
\setcounter{section}{1}

\begin{center}
{\bf 1. \  Introduction}
\end{center}
\renewcommand{\theequation}{\arabic{section}.\arabic{equation}}

The double hypergeometric series, defined by Humbert \cite{Manako1,Manako2} are the following :
\begin{equation} \label{1.1}
\Psi_2(a;\,b,\,c;\,x,\,y)=\sum_{n,k=0}^\infty\frac{(a)_{n+k}}{(b)_{n}\;(c)_{k}}\frac{x^n\;y^k}{n!\;k!}, 
\end{equation}
\begin{equation} \label{1.2}
\Phi_3(b;\,c;\,x,\,y)=\sum_{n,k=0}^\infty\frac{(b)_n}{(c)_{n+k}}\frac{x^n\;y^k}{n!\;k!}
\end{equation}
The double series \eqref{1.1} and \eqref{1.2} converge absolutely for all $x,\,y \in \mathbb{C}$.
Very recently, Manako \cite{Manako3} established a few results for the series $\Phi_2$, $\Phi_3$ and $\Psi_2$ out of which, two results are given here :

For  $b$, $c \neq 0, -1, -2,...$
\begin{equation}
\label{1.3}
\Psi_2(b;\,b,\,c;\,x,\,y)=\exp(x+y)\,\Phi_2(c-b;\;c;\;-y,xy)
\end{equation}
and  for $b$, $c \neq 0, -1, -2,...$ and $|x|\neq 0$
\begin{equation}
\label{1.4}
\Psi_2(a;\,b,\,c;\,x,\;y) = \sum_{k=0}^\infty\frac{(a)_k}{(b)_k}\, {}_2F_1\left[ \begin{array}{c} -k,\,-k-b+1 \\  c \end{array} ;\,\frac{y}{x} \right]\,\frac{x^k}{k!}\,.
\end{equation}
In the same paper, using \eqref{1.3} and \eqref{1.4}, Manako \cite{Manako3} established the following new result for $\Phi_3$ in terms of series of $_2F_1$,
\begin{equation}
\label{1.5}
\Phi_3(b;\,c;\,x,\,y)= \exp \left( x+\frac{y}{x} \right)\,\sum_{k=0}^\infty\,\frac{1}{k!}\, \left(-\frac{y}{x}\right)^k 
 {}_2F_1\left[ \begin{array}{c}  -k,\,-k-b+1 \\  c  \end{array} ;\,\frac{x^2}{y} \right] \, .
\end{equation}
In 2013, Rathie[4] obtained the following result for the series $\Phi_2$ : \\
For $ c \neq 0, -1, -2, \cdots$
\begin{equation}         
\Phi_2(a, b; c; x, y) = \sum_{m=0}^\infty \frac{(a)_m}{(c)_m} \, {}_2F_1\left[ \begin{array}{c} -m,\,b   \\ 1-a-m  \end{array} ;\,\frac{y}{x} \right] \, \frac{x^m}{m!} .
\end{equation}
and discussed some special cases.\\

The aim of this short research note is to establish \eqref{1.5} by another method. In the end, we consider some interesting special cases.
\section{Derivation of (1.5)}
In order to derive \eqref{1.5}, let us denote its right-hand side by $S$ and expressing $_2F_1$ as a series, we have
$$S=\exp \left( x+\frac{y}{x} \right)\, \sum_{k=0}^\infty\,\sum_{m=0}^k\,\frac{(-k)_m\,(-k-c+b+1)_m\,(-1)^k}{(c)_m\, m!\,k!}\,x^{2m-k}\,y^{k-m}$$
Using
$$\frac{(-k)_m}{k!}=\frac{(-1)^m}{(k-m)!}$$
we have
$$S=\exp \left( x+\frac{y}{x} \right)\,\sum_{k=0}^\infty\,\sum_{m=0}^k\,\frac{(-k-c+b+1)_m}{(c)_m\,m!}\,\frac{(-1)^{k+m}\, x^{2m-k}\,y^{k-m}}{(k-m)!}$$
Replacing $k$ by $k+m$, and using the result[5, Lemma 10, equ. 2, p-57]
$$S= \exp\left( x+\frac{y}{x} \right)\,\sum_{k=0}^\infty\,\sum_{m=0}^\infty\,\frac{(-k-m-c+b+1)_m}{\,(c)_m\,m!}\,\frac{(-1)^{k}\,x^{m-k}\,y^k}{k!}$$
Using $(-k-m-c+b+1)_m=(-1)^m\;(k+c-b)_m$, we have
$$S=\exp\left( x+\frac{y}{x} \right)\,\sum_{k=0}^\infty\,\sum_{m=0}^\infty\,\frac{(-1)^m\;(k+c-b)_m}{\,(c)_m\;m!}\,\frac{(-1)^{k}\,x^{m-k}\;y^k}{k!}$$
which can be written as
$$S=\exp \left( x+\frac{y}{x} \right)\,\sum_{k=0}^\infty\,\frac{(-1)^{k}}{k!}\, \left(\frac{y}{x} \right)^k\,
\sum_{m=0}^\infty\,\frac{(-1)^m\;(k+c-b)_m}{\,(c)_m\;m!}\,x^m $$
Now, summing up the inner series, we have
\begin{align*}
&S= \exp \left( x+\frac{y}{x} \right) \,\sum_{k=0}^\infty\,\frac{(-1)^{k}}{k!}\,\left(\frac{y}{x} \right)^k 
{}_1F_1\left[ \begin{array}{c}  k+c-b \\  c \end{array} ;\,-x \right] \, .\\
&\,\,\,=\exp \left(\frac{y}{x} \right)\,\sum_{k=0}^\infty\,\frac{\left(-\frac{y}{x} \right)^k}{k!}\,\,\left\{\exp (x)\,
{}_1F_1\left[ \begin{array}{c} k+c-b \\ c  \end{array} ;\,-x \right] \right\} \, .
\end{align*}
Using Kummer's first transformation \cite{Manako4}
\begin{equation}
\label{2.6}
\exp (x)\,\,{}_1F_1\left[\begin{array}{c}  c-b \\ c  \end{array} ;\, -x \right] \, = 
{}_1F_1\left[\begin{array}{c}  b \\ c  \end{array} ;\,x \right] \, .
\end{equation}
we have 
$$S=\exp\left(\frac{y}{x}\right)\,\sum_{k=0}^\infty\,\frac{\left(-\frac{y}{x}\right)^k}{k!}
{}_1F_1\left[ \begin{array}{c}  b-k \\ c  \end{array} ;\,x \right] \, .
$$
Again, expressing $_1F_1$ as a series, we have
$$S=\exp\left(\frac{y}{x}\right)\,\sum_{k=0}^\infty\,\sum_{n=0}^\infty\,\frac{\left(-\frac{y}{x}\right)^k}{k!}\;\frac{(b-k)_n}{(c)_n} \frac{x^n}{n!}$$
Since $(b-k)_k=\frac{(b)_n\,(1-b)_n}{(1-b-n)_k}$, we have, therefore
$$S=\exp\left(\frac{y}{x}\right)\,\sum_{n=0}^\infty\,\frac{(b)_n}{(c)_n}\;\frac{x^n}{n!}\,\sum_{k=0}^\infty\,\frac{(1-b)_k}{(1-b-n)_k}\;\frac{\left(-\frac{y}{x}\right)^k}{k!}$$
Summing up the inner series
\begin{align*}
& S=\exp\left(\frac{y}{x}\right)\;\sum_{n=0}^\infty\,\frac{(b)_n}{(c)_n}\;\frac{x^n}{n!}\,
{}_1F_1\left[ \begin{array}{c} 1-b   \\  1-b-n \end{array} ;\, -\frac{y}{x} \right] \\
&=\sum_{n=0}^\infty\,\frac{(b)_n}{(c)_n}\,\frac{x^n}{n!}\,\left\{\exp\left(\frac{y}{x}\right)\,
{}_1F_1\left[ \begin{array}{c} 1-b \\ 1-b-n  \end{array} ;\,-\frac{y}{x} \right] \right\} \,.
\end{align*}
Using \eqref{2.6}, we have
$$S=\sum_{n=0}^\infty \,\frac{(b)_n}{(c)_n}\;\frac{x^n}{n!}\,
{}_1F_1\left[ \begin{array}{c}  -n \\ 1-b-n  \end{array} ;\,\frac{y}{x} \right] \,.
$$
Expressing $_1F_1$ as a series
$$ S=\sum_{n=0}^\infty\,\sum_{m=0}^n\,\frac{(b)_n}{(c)_n}\;\frac{(-n)_m}{(1-b-n)_m}\;\frac{x^n\,\left(\frac{y}{x}\right)^m}{m!\;n!}$$
Using $(-n)_m=\frac{(-1)^m\;n!}{(n-m)!}$
$$ S=\sum_{n=0}^\infty\,\sum_{m=0}^n\,\frac{(b)_n\;(-1)^m}{(c)_n\;(1-b-n)_m}\,\frac{x^n\;\left(\frac{y}{x}\right)^m}{m!\;(n-m)!}$$
Changing $n$ to $n+m$ and using [5, Lemma 10, equ. 2, p-57 ], we have
$$ S=\sum_{n=0}^\infty\,\sum_{m=0}^\infty\,\frac{(b)_{n+m}\;(-1)^m\;y^m\;x^n}{(c)_{n+m}\;(1-b-n-m)_m\;m!\;n!}$$
Using $(1-b-n-m)_m= \frac{(-1)^m\;(b)_{n+m}}{(b)_n}$
$$ S=\sum_{n=0}^\infty\,\sum_{m=0}^\infty\,\frac{(b)_n\;x^n\;y^m}{(c)_{n+m}\;m!\;n!}$$
Finally, using definition \eqref{1.2}, we have
$$ S=\Phi_3(b;\; c;\; x,\; y) $$
This completes the proof of \eqref{1.5}.
\section{SPECIAL CASES}

In this section, we shall mention two interesting special cases of our results \eqref{1.5}.

In \eqref{1.5}, if we take $y=x^2$, we have
\begin{equation}
\label{3.1}
\Phi_3(b;\,c;\,x,\,x^2)=\exp\left(2x \right)\; \sum_{k=0}^\infty\;\frac{(-x)^k}{k!}\; \left(-\frac{y}{x} \right)^k
{}_2F_1\left[ \begin{array}{c}  -k,\; -k-c+b+1  \\ c   \end{array} ; \; 1 \right]\,.
\end{equation}
The $_2F_1$ appearing on the right-hand side of \eqref{3.1} can be evaluated with the help of classical Gauss's summation theorem \cite{Manako5}
\begin{equation}
\label{3.2}
_2F_1 \left[\begin{array}{c} a, \, b \\ c \end{array} ; 1 \right]=\frac{\Gamma(c)\;\Gamma(c-a-b)}{\Gamma(c-a)\;\Gamma(c-b)}
\end{equation}
provided $\text{Re}(c-a-b)>0$.\\
and we get, after some simplification, the following new reduction formula
\begin{equation}
\label{3.3}
\Phi_3(b;\,c;\,x,\,x^2)=\exp\left(2x\right)\,{}_2F_2 \left[ \begin{array}{c}  c-\frac{b}{2}, \, c - \frac{b}{2}-\frac{1}{2}\\  \,c,\,2c-b-1\end{array} \,;\, -4x\,\right]
\end{equation}
Further, in \eqref{3.3}, if we take $c=2b$, we have
$$\Phi_3(b;\,2b;\,x,\,x^2)=\exp\left(2x\right)\,{}_2F_2 \left[\begin{array}{c} \frac{3b}{2}, \; \frac{3b-1}{2} \\ 2b,\; 3b-1 \end{array};\;-4x\,\right] $$
and using \eqref{1.4} after simplification, we get
\begin{equation}
\label{3.4}
\Psi_2(b;\; b,\;2b;\;x,\; x)=\,_2F_2 \left[\begin{array}{c} \frac{3b}{2},\; \frac{3b-1}{2} \\ 2b, \; 3b-1\end{array};\; 4x\right]
\end{equation}
which is a special case of the following result
$$\Psi_2(b;\, b,\, c;\, x,\, x)=\,_3F_3 \left[\begin{array}{c} a,\; \frac{c+b}{2},\; \frac{c+b-1}{2} \\  b,\;c,\; c+b-1\end{array} ;\;4x\right]$$
given by Burchnall and Chaundy \cite{Manako5,Manako6}, also recorded in \cite{Manako1,Manako2}.

Similarly, other results can also be obtained. 

\section*{Remark :}
In 2015, Choi and Rathie [8] obtained explicit expressions of 

$\Phi_2( a,\; a+i;\; c;\; x, \; -x)\,$

and

$\,\Psi_2(a, \; c;\; c+i; \; x,\;  -x)\,$

each for $\,i= 0, \pm 1, \pm 2, \cdots, \pm 5\,$.\\
and deduced interesting summation formulas.


\vspace{0.6cm}

\begin{thebibliography}{99}
\footnotesize{
\bibitem {Manako1}
Brychkov Yu.\,A., {\em Handbook of Special Functions: Derivatives, Integrals, Series and Other Formulas}, CRC Press, Boca Raton, Fl,  2008.
\bibitem {Manako2}
Prudnikov A.\,P., Brychkov Yu.\,A., and Marichev O.\,I.,  {\em Table of Integrals, Series and Products}, Academic Press, New York, 2007.
\bibitem {Manako3}
Manako V.\,V., {\em A connection formula between double hypergeometric series $\Psi_2$ and $\Phi_3$}, Integral Transforms and Special Functions, 23(7), 503-508, 2012.
\bibitem{Rathie} Rathie, A.K., {\em On a representation of Humberts double hypergeometric series $\Phi_2$ in a series of Gauss's $_{2}F_1$ function}, arxiv : 1312.0064v1[math.CV], 30 Nov. 2013.
\bibitem {Manako4}
Rainville E.\,D., {\em Special Functions}, The Macmillian company, New York, 1960.
\bibitem {Manako5}
Burchnall J.\,L., Chaundy T.\,W., {\em Expansions of Appell's double hypergeometric functions-I}, Quarterly J. Math., Oxford, 11, 249-270, 1940.
\bibitem {Manako6}
Burchnall J.\,L., Chaundy T.\,W. , {\em Expansions of Appell's double hypergeometric functions-II}, Quarterly J. Math., Oxford, 12, 112-128, 1941.
\bibitem{Choi} Choi, J and Rathie, A.K., {\em Certain summation formulas for Humbert's double hypergeometric series $\Psi_2$ and $\Phi_2$}, Commun. Korean Math. Soc. 30(4), 439-446, (2015).
}
\end{thebibliography}
\end{document}